\input amstex
\magnification\magstephalf
\documentstyle{amsppt}

\hsize 5.72 truein
\vsize 7.9 truein
\hoffset .39 truein
\voffset .26 truein
\mathsurround 1.67pt
\parindent 20pt
\normalbaselineskip 13.8truept
\normalbaselines
\binoppenalty 10000
\relpenalty 10000
\csname nologo\endcsname 


\font\bc=cmb10
\font\tenbsy=cmbsy10

\catcode`\@=11

\def\myitem#1.{\item"(#1)."\advance\leftskip10pt\ignorespaces}

\def\qedsymbol{{\mathsurround\z@$\square$}}
\redefine\qed{\relaxnext@\ifmmode\let\next\@qed\else
  {\unskip\nobreak\hfil\penalty50\hskip2em\null\nobreak\hfil
    \qedsymbol\parfillskip\z@\finalhyphendemerits0\par}\fi\next}
\def\@qed#1$${\belowdisplayskip\z@\belowdisplayshortskip\z@
  \postdisplaypenalty\@M\relax#1
  $$\par{\lineskip\z@\baselineskip\z@\vbox to\z@{\vss\noindent\qed}}}
\outer\redefine\beginsection#1#2\par{\par\penalty-250\bigskip\vskip\parskip
  \leftline{\tenbsy x\bf#1. #2}\nobreak\smallskip\noindent}
\outer\redefine\genbeginsect#1\par{\par\penalty-250\bigskip\vskip\parskip
  \leftline{\bf#1}\nobreak\smallskip\noindent}

\def\next{\let\@sptoken= }\def\next@{ }\expandafter\next\next@
\def\@futureletnext#1{\let\nextii@#1\futurelet\next\@flti}
\def\@flti{\ifx\next\@sptoken\let\next@\@fltii\else\let\next@\nextii@\fi\next@}
\expandafter\def\expandafter\@fltii\next@{\futurelet\next\@flti}

\let\zeroindent\z@
\let\savedef@\endproclaim\let\endproclaim\relax 
\define\chkproclaim@{\add@missing\endroster\add@missing\enddefinition
  \add@missing\endproclaim
  \envir@stack\endproclaim
  \edef\endit@{\leftskip\the\leftskip\rightskip\the\rightskip}}
\let\endproclaim\savedef@
\def\thing@{.\enspace\egroup\ignorespaces}
\def\thingi@(#1){ \rm(#1)\thing@}
\def\thingii@\cite#1{ \rm\@pcite{#1}\thing@}
\def\thingiii@{\ifx\next(\let\next\thingi@
  \else\ifx\next\cite\let\next\thingii@\else\let\next\thing@\fi\fi\next}
\def\thing#1#2#3{\chkproclaim@
  \ifvmode \medbreak \else \par\nobreak\smallskip \fi
  \noindent\advance\leftskip#1
  \hskip-#1#3\bgroup\bc#2\unskip\@futureletnext\thingiii@}
\let\savedef@\endproclaim\let\endproclaim\relax 
\def\endit{\endproclaim\endit@\let\endit@\undefined}
\let\endproclaim\savedef@

\def\lemma#1{\thing\parindent{Lemma #1}\sl}
\def\prop#1{\thing\parindent{Proposition #1}\sl}
\def\thm#1{\thing\parindent{Theorem #1}\sl}

\def\narrowthing#1{\chkproclaim@\medbreak\narrower\noindent
  \it\def\next{#1}\def\next@{}\ifx\next\next@\ignorespaces
  \else\bgroup\bc#1\unskip\let\next\narrowthing@\fi\next}
\def\narrowthing@{\@futureletnext\thingiii@}

\def\@cite#1,#2\end@{{\rm([\bf#1\rm],#2)}}
\def\cite#1{\in@,{#1}\ifin@\def\next{\@cite#1\end@}\else
  \relaxnext@{\rm[\bf#1\rm]}\fi\next}
\def\@pcite#1{\in@,{#1}\ifin@\def\next{\@cite#1\end@}\else
  \relaxnext@{\rm([\bf#1\rm])}\fi\next}

\advance\minaw@ 1.2\ex@
\atdef@[#1]{\ampersand@\let\@hook0\let\@twohead0\brack@i#1,\z@,}
\def\brack@{\z@}
\let\@@hook\brack@
\let\@@twohead\brack@
\def\brack@i#1,{\def\next{#1}\ifx\next\brack@
  \let\next\brack@ii
  \else \expandafter\ifx\csname @@#1\endcsname\brack@
    \expandafter\let\csname @#1\endcsname1\let\next\brack@i
    \else \Err@{Unrecognized option in @[}%
  \fi\fi\next}
\def\brack@ii{\futurelet\next\brack@iii}
\def\brack@iii{\ifx\next>\let\next\brack@gtr
  \else\ifx\next<\let\next\brack@less
    \else\relaxnext@\Err@{Only < or > may be used here}
  \fi\fi\next}
\def\brack@gtr>#1>#2>{\setboxz@h{$\m@th\ssize\;{#1}\;\;$}%
 \setbox@ne\hbox{$\m@th\ssize\;{#2}\;\;$}\setbox\tw@\hbox{$\m@th#2$}%
 \ifCD@\global\bigaw@\minCDaw@\else\global\bigaw@\minaw@\fi
 \ifdim\wdz@>\bigaw@\global\bigaw@\wdz@\fi
 \ifdim\wd@ne>\bigaw@\global\bigaw@\wd@ne\fi
 \ifCD@\enskip\fi
 \mathrel{\mathop{\hbox to\bigaw@{$\ifx\@hook1\lhook\mathrel{\mkern-9mu}\fi
  \setboxz@h{$\displaystyle-\m@th$}\ht\z@\z@
  \displaystyle\m@th\copy\z@\mkern-6mu\cleaders
  \hbox{$\displaystyle\mkern-2mu\box\z@\mkern-2mu$}\hfill
  \mkern-6mu\mathord\ifx\@twohead1\twoheadrightarrow\else\rightarrow\fi$}}%
 \ifdim\wd\tw@>\z@\limits^{#1}_{#2}\else\limits^{#1}\fi}%
 \ifCD@\enskip\fi\ampersand@}
\def\brack@less<#1<#2<{\setboxz@h{$\m@th\ssize\;\;{#1}\;$}%
 \setbox@ne\hbox{$\m@th\ssize\;\;{#2}\;$}\setbox\tw@\hbox{$\m@th#2$}%
 \ifCD@\global\bigaw@\minCDaw@\else\global\bigaw@\minaw@\fi
 \ifdim\wdz@>\bigaw@\global\bigaw@\wdz@\fi
 \ifdim\wd@ne>\bigaw@\global\bigaw@\wd@ne\fi
 \ifCD@\enskip\fi
 \mathrel{\mathop{\hbox to\bigaw@{$%
  \setboxz@h{$\displaystyle-\m@th$}\ht\z@\z@
  \displaystyle\m@th\mathord\ifx\@twohead1\twoheadleftarrow\else\leftarrow\fi
  \mkern-6mu\cleaders
  \hbox{$\displaystyle\mkern-2mu\copy\z@\mkern-2mu$}\hfill
  \mkern-6mu\box\z@\ifx\@hook1\mkern-9mu\rhook\fi$}}%
 \ifdim\wd\tw@>\z@\limits^{#1}_{#2}\else\limits^{#1}\fi}%
 \ifCD@\enskip\fi\ampersand@}

\def\pr@m@s{\ifx'\next\let\nxt\pr@@@s \else\ifx^\next\let\nxt\pr@@@t
  \else\let\nxt\egroup\fi\fi \nxt}

\define\widebar#1{\mathchoice
  {\setbox0\hbox{\mathsurround\z@$\displaystyle{#1}$}\dimen@.1\wd\z@
    \ifdim\wd\z@<.4em\relax \dimen@ -.16em\advance\dimen@.5\wd\z@ \fi
    \ifdim\wd\z@>2.5em\relax \dimen@.25em\relax \fi
    \kern\dimen@ \overline{\kern-\dimen@ \box0\kern-\dimen@}\kern\dimen@}%
  {\setbox0\hbox{\mathsurround\z@$\textstyle{#1}$}\dimen@.1\wd\z@
    \ifdim\wd\z@<.4em\relax \dimen@ -.16em\advance\dimen@.5\wd\z@ \fi
    \ifdim\wd\z@>2.5em\relax \dimen@.25em\relax \fi
    \kern\dimen@ \overline{\kern-\dimen@ \box0\kern-\dimen@}\kern\dimen@}%
  {\setbox0\hbox{\mathsurround\z@$\scriptstyle{#1}$}\dimen@.1\wd\z@
    \ifdim\wd\z@<.28em\relax \dimen@ -.112em\advance\dimen@.5\wd\z@ \fi
    \ifdim\wd\z@>1.75em\relax \dimen@.175em\relax \fi
    \kern\dimen@ \overline{\kern-\dimen@ \box0\kern-\dimen@}\kern\dimen@}%
  {\setbox0\hbox{\mathsurround\z@$\scriptscriptstyle{#1}$}\dimen@.1\wd\z@
    \ifdim\wd\z@<.2em\relax \dimen@ -.08em\advance\dimen@.5\wd\z@ \fi
    \ifdim\wd\z@>1.25em\relax \dimen@.125em\relax \fi
    \kern\dimen@ \overline{\kern-\dimen@ \box0\kern-\dimen@}\kern\dimen@}%
  }

\catcode`\@\active

\let\PVstyle=d 

\font\tenscr=rsfs10 
\font\sevenscr=rsfs7 
\font\fivescr=rsfs5 
\skewchar\tenscr='177 \skewchar\sevenscr='177 \skewchar\fivescr='177
\newfam\scrfam \textfont\scrfam=\tenscr \scriptfont\scrfam=\sevenscr
\scriptscriptfont\scrfam=\fivescr
\define\scr#1{{\fam\scrfam#1}}
\let\Cal\scr

\let\0\relax 
\define\restrictedto#1{\big|_{#1}}

\loadbold

\define\exc{{\text{exc}}}
\define\Gm{\Bbb G_{\text m}}
\define\Spec{\operatorname{Spec}}
\define\Supp{\operatorname{Supp}}
\mathchardef\idot="202E

\outer\redefine\subhead#1\par{\par\medbreak\leftline{\bc #1}\smallskip\noindent}

\topmatter
\title Transplanting Faltings' Garden\endtitle
\author Paul Vojta\endauthor
\affil University of California, Berkeley\endaffil
\address Department of Mathematics, University of California,
  970 Evans Hall\quad\#3840, Berkeley, CA \ 94720-3840\endaddress
\date 14 January 2009 \enddate
\thanks This paper was written while the author enjoyed the kind hospitality
of the Fields Institute; partially supported by NSF grant DMS-0753152.\endthanks
\subjclassyear{2000}

\abstract In his contribution to the {\it Baker's Garden\/} book, Faltings
gives a family of examples of irreducible divisors $D$ on $\Bbb P^2$ for which
$\Bbb P^2\setminus D$ has only finitely many integral points over any given
localization of a number ring away from finitely many places.  He also notes
that neither $\Bbb P^2\setminus D$ nor the \'etale covers used in his proof
embed into semiabelian varieties, so his examples do not easily reduce to
known results about such subvarieties.  In this note, we show how Faltings'
results follow directly from a theorem of Evertse and Ferretti; hence these
examples can be explained by noting that if one pulls back to a cover of
$\Bbb P^2$ \'etale outside of $D$ and then adds components to the pull-back
of $D$ then one can embed the complement into a semiabelian variety and obtain
useful diophantine approximation results for the original divisor $D$.
\endabstract
\endtopmatter

\document

In his contribution to the {\it Baker's Garden\/} volume, Faltings \cite{F}
gave examples of irreducible divisors $D$ on $\Bbb P^2$ for which
$\Bbb P^2\setminus D$ has only finitely many integral points over any
number ring, and over any localization of such a ring away from finitely many
places.  This was further explored by Zannier \cite{Z} using methods of
Zannier and Corvaja, although he uses a different family of examples.
This family has substantial overlap with the examples of Faltings but
does not contain all of his examples.

The present note explores Faltings' examples from the point of view of work
of Evertse and Ferretti \cite{E-F}.  The main theorem, Theorem \02.6,
covers all of Faltings' examples, yet its proof follows rather directly from
the main theorem of Evertse and Ferretti.  Since the latter theorem relies
on Schmidt's Subspace Theorem, it necessarily involves varieties that can
be embedded into semiabelian varieties (actually $\Gm^N$).
Thus, Faltings' examples can be viewed as examples where one adds components
to the divisor $D$ to obtain a divisor whose complement can be embedded
into a semiabelian variety, and then the resulting diophantine inequality
implies a useful inequality for the original divisor $D$.

This is an approach that deserves further attention.

In particular, the Shafarevich conjecture (on semistable abelian varieties
over a given number field with good reduction outside of a given finite set
of places, proved by Faltings in 1983) stands out as presently the only
diophantine result with all the hallmarks of a result proved by Thue's method
(ineffective, but with bounds on the number of counterexamples), but which
has not been proved by Thue's method.  The relevant variety, $\Cal A_{g,n}$,
also cannot be embedded into a semiabelian variety.  It would be interesting
to know if a proof via Thue's method could be obtained by adding divisor
components on some \'etale cover to obtain an embedding into a semiabelian
variety.

The first section of this paper briefly summarizes the geometric setting
in Faltings' paper.  Section \02 states and proves the main theorem of this
paper.  Finally, the third section briefly gives the counterpart to the
main theorem in the case of holomorphic curves.

\beginsection{\01}{The Geometric Setting}

This section describes the geometric setting used in Falting's paper \cite{F}.
This setting is only sketched; for full details see his paper.

Faltings starts with a smooth geometrically irreducible algebraic surface
over a field $k$ of characteristic zero, and an ample line sheaf $\Cal L$
on $X$.  Under certain assumptions on $\Cal L$ (satisfied if $\Cal L$ is a
tensor product of five very ample line sheaves and if
$\Cal K_X\otimes\Cal L^{\otimes3}$ is ample), generic three-dimensional
subspaces of $\Gamma(X,\Cal L)$ (as determined by a dense open subset
of the appropriate Grassmannian) determine finite morphisms $f\:X\to\Bbb P^2$
that satisfy the following conditions.
\roster
\myitem i.  The ramification locus $Z$ of $f$ is smooth and irreducible,
and the ramification index is $2$.
\myitem ii.  The restriction of $f$ to $Z$ is birational onto its
image $D\subseteq\Bbb P^2$.
\myitem iii.  $D$ is nonsingular except for cusps and simple double points.
\myitem iv.  Let $Y\to X\to\Bbb P^2$ denote the Galois closure of
$X\to\Bbb P^2$, and let $n=\deg f$.  Then $Y$ is
smooth and the Galois group is the full symmetric group $\Cal S_n$.
\myitem v.  The ramification locus of $Y$ over $\Bbb P^2$ is the sum of
distinct conjugate effective divisors $Z_{ij}$, $1\le i<j\le n$.  They have
smooth supports, and are disjoint with the following two exceptions.
Points of $Y$ lying over double points of $D$ are fixed points of a subgroup
$\Cal S_2\times\Cal S_2$ of $\Cal S_n$, and they lie on $Z_{ij}\cap Z_{\ell m}$
with distinct indices $i,j,\ell,m$.  Points of $Y$ lying over cusps of $D$
are fixed points of a subgroup $\Cal S_3$ of $\Cal S_n$, and lie on
$Z_{ij}\cap Z_{i\ell}\cap Z_{j\ell}$.
\endroster

For convenience, write $Z_{ij}=Z_{ji}$ when $i,j\in\{1,\dots,n\}$ and $i>j$,
and let
$$A_i = \sum_{j\ne i} Z_{ij}\;.$$
Let $L$ denote the divisor class of $\Cal L$ on $X$, identified with a divisor
class on $Y$ via pull-back.  We then have
$$2\sum_{i<j}Z_{ij} = \sum A_i \sim dL\tag\01.1$$
on $Y$.

Under the additional assumptions that $dL-\alpha A_i$ is ample for some
rational $\alpha>12$ and all $i$, Faltings then shows that there are only
finitely many integral points on $\Bbb P^2\setminus D$.  (In his proof, the
only explicit lower bound on $\alpha$ occurs at the end of Section 3 (page 242),
where $\alpha>6$ is used.  However, by (\01.1), the bound on the index
just prior to that point should really be that $4rd\cdot h(y_1)$ is bounded
by $\alpha rd/3\cdot h(y_1)+\text{constant}$, so one one really needs
$\alpha>12$ there.)

Faltings' proof proceeds by first noting that $Y\setminus\bigcup Z_{ij}$ is
an \'etale cover of $\Bbb P^2\setminus D$, so by \cite{S, \S\,4.2} there is
a fixed number field $k'$ such that integral points on $\Bbb P^2\setminus D$
over $k$ lift to integral points on $Y\setminus\bigcup Z_{ij}$ over $k'$.
He then shows finiteness of integral points on the latter, by using his
Product Theorem together with some estimates of expectation values.

Finally, Faltings notes that if $X=\Bbb P^1\times\Bbb P^1$ and
$\Cal L=\Cal O(a,b)$ with $a$ and $b$ coprime and $a,b\ge5$, then
the above conditions are satisfied, and moreover $Y\setminus\bigcup Z_{ij}$
does not embed into a semiabelian variety.

\beginsection{\02}{The Theorem}

This section states and proves the theorem that can be obtained by applying
the theorem of Evertse and Ferretti.

For the purposes of this section, a {\bc variety} over a field $k$ is
an integral scheme, separated and of finite type over $\Spec k$.  If $k$
is a number field and $v$ is a place of $k$, then $\|x\|_v$ for $x\in k$
is as defined in \cite{V~1} or \cite{V~2}.

First we recall the theorem of Evertse and Ferretti \cite{E-F, Theorem 1.1}:

\thm{\02.1}  Let $k$ be a number field, let $S$ be a finite set of places
of $k$, let $X$ be a closed subvariety of $\Bbb P^N_k$ of dimension $n\ge1$,
and let $0<\epsilon\le1$.  Further, for each $v\in S$
let $f_0^{(v)},\dots,f_n^{(v)}$ be a system of homogeneous polynomials
in $k[x_0,\dots,x_N]$ satisfying
$$X\cap Z(f_0^{(v)})\cap\dots\cap Z(f_n^{(v)})=\emptyset\tag\02.1.1$$
and $Z(f_i^{(v)})\nsupseteq X$ for all $i$.  Then the inequality
$$-\sum_{v\in S}\sum_{i=0}^n
    \log\frac{\|f_i^{(v)}(\bold x)\|_v^{1/\deg f_i^{(v)}}}{\|\bold x\|_v}
  \le (n+1+\epsilon)h_k(\bold x)\tag\02.1.2$$
holds for all $\bold x\in X(k)$ outside of a proper Zariski-closed subset
of $X$.  Here
$$\|\bold x\|_v = \max\{\|x_0\|_v,\dots,\|x_N\|_v\}$$
where $[x_0:\dots,x_N]$ are homogeneous coordinates in $k$ for $\bold x$.
This of course depends on the choice of homogeneous coordinates, but so does
$\|f_i^{(v)}(\bold x)\|_v$, so the fractions in (\02.1.2) are well defined.
\endit

In the present application, the polynomials $f_i^{(v)}$ will all be
linear polynomials, associated to hyperplanes $H_i^{(v)}$ in $\Bbb P^N_k$,
so the left-hand side of (\02.1.2) can be expressed in terms of Weil functions
$\lambda_{H_i^{(v)},v}$:
$$\sum_{v\in S} \sum_{i=0}^n \lambda_{H_i^{(v)},v}(\bold x)
  \le (n+1+\epsilon)h_k(\bold x)\;.\tag\02.2$$
More generally, let $H_1,\dots,H_q$ be hyperplanes in $\Bbb P^N_k$,
and let $\Cal J$ be the set of $(n+1)$\snug-element subsets $J$
of $\{1,\dots,q\}$ for which
$$X\cap \bigcap_{j\in J}H_j=\emptyset\;.$$
If $\Cal J$ is not empty, then (\02.2) can be restated as
$$\sum_{v\in S} \max_{J\in\Cal J} \sum_{j\in J} \lambda_{H_j,v}(\bold x)
  \le (n+1+\epsilon)h_k(\bold x)\;.\tag\02.3$$
Indeed, at each place there are only finitely many choices for $J$, so this
follows by invoking (\02.2) finitely many times.

The bulk of this section is devoted to finding hyperplanes in a suitable
projective embedding that allow one to relate Weil functions relative to
the $Z_{ij}$ to the left-hand side of (\02.3).

The geometric setting under consideration can be summarized as follows.

\narrowthing{\02.4}  Let $Y$ be a projective surface over a field $k$,
let $n$ be a positive integer, let $Z_{ij}$ be effective Cartier divisors
on $Y$ for all $1\le i<j\le n$, and let $\alpha$ be a rational number.
Assume that the supports of the $Z_{ij}$ are disjoint, except that $Z_{ij}$
may meet $Z_{\ell m}$ in finitely many points if $\{i,j\}$ and $\{\ell,m\}$
are disjoint, and that if $1\le i<j<\ell\le n$, then
$$Z_{ij}\cap Z_{i\ell}\cap Z_{j\ell} = Z_{ij}\cap Z_{i\ell}
  = Z_{ij}\cap Z_{j\ell} = Z_{i\ell}\cap Z_{j\ell}\;,$$
and this set is finite (possibly empty).
For $1\le j<i\le n$ write $Z_{ij}=Z_{ji}$, let
$$A_i = \sum_{j\ne i} Z_{ij}$$
for all $i$, and let
$$M = \sum A_i = 2\sum_{i<j} Z_{ij}\;.$$
Finally, assume that $M$ is ample and that $M-\alpha A_i$ is an ample
$\Bbb Q$\snug-divisor for all $i$.
\endit

The bulk of the work in this section consists of proving the following result.

\prop{\02.5}  Let $k$ be a local or global field, and let $Y$, $n$,
$\{Z_{ij}\}_{i<j}$, $\{A_i\}_i$, $M$, and $\alpha$ be as in (\02.4).
Assume that $n\ge4$.  Let $v$ be a place of $k$, and fix Weil functions
$\lambda_{ij,v}$ for each $Z_{ij}$ at $v$.  Let $\beta$ be an integer such that
$\beta\alpha\in\Bbb Z$ and such that $\beta M$ and all $\beta(M-\alpha A_i)$
are very ample.  Fix an embedding $Y\hookrightarrow\Bbb P^N_k$ associated
to a complete linear system of $\beta M$, and regard $Y$ as a subvariety
of $\Bbb P^N_k$ via this embedding.  Then
\roster
\myitem a.  There is a finite list $H_1,\dots,H_q$ of hyperplanes in
$\Bbb P^N_k$, with associated Weil functions $\lambda_{H_j,v}$ at $v$
for all $j$, with the following property.  Let $\Cal J$ be the set of
$3$\snug-element subsets $J=\{j_0,j_1,j_2\}$ of $\{1,\dots,q\}$ for which
$Y\cap H_{j_0}\cap H_{j_1}\cap H_{j_2}=\emptyset$.  Then $\Cal J\ne\emptyset$,
and the inequality
$$\max_{J\in\Cal J}\sum_{j\in J} \lambda_{H_j,v}(y)
  \ge \beta\alpha \sum_{i<j}\lambda_{ij,v}(y) + O(1)\tag\02.5.1$$
holds for all $y\in Y(k)$ not lying on the support of any $Z_{ij}$ or
on any of the $H_j$.
\myitem b.  Let $C$ be an integral curve in $Y$, not contained in the support
of any $Z_{ij}$.  Then there is a finite list $H_1,\dots,H_q$ of hyperplanes,
with associated Weil functions as before, with the following property.
Let $\Cal J$ be the set of all $2$\snug-element subsets $J=\{j_0,j_1\}$
of $\{1,\dots,q\}$ for which $C\cap H_{j_0}\cap H_{j_1}=\emptyset$.
Then $\Cal J\ne\emptyset$, and the inequality
$$\max_{J\in\Cal J}\sum_{j\in J} \lambda_{H_j,v}(y)
  \ge \frac{\beta\alpha}{2} \sum_{i<j}\lambda_{ij,v}(y) + O(1)\tag\02.5.2$$
holds for all but finitely many $y\in C(k)$.
\endroster
In each case the implicit constant in $O(1)$ is independent of $y$ but may
depend on all of the other data.
\endit

\demo{Proof}  The proof relies mainly on two lemmas.  These lemmas replace
Faltings' computations of ideals associated to indices.

\lemma{\02.5.3}  Let $i,j,\ell,m$ be distinct indices.  Then:
\roster
\myitem a.  there exist hyperplanes $H_0$, $H_1$, and $H_2$ in $\Bbb P^N_k$,
such that
$$Y\cap H_0\cap H_1\cap H_2=\emptyset$$
and
$$\lambda_{H_0,v}(y) + \lambda_{H_1,v}(y) + \lambda_{H_2,v}(y)
  \ge \beta\alpha(\lambda_{ij,v}(y) + \lambda_{\ell m,v}(y)) + O(1)
  \tag\02.5.3.1$$
for all $v\in S$ and all $y\in Y(k)$ outside of $H_0\cup H_1\cup H_2$; and
\myitem b.  given any integral curve $C\subseteq Y$ not contained in any of the
$Z_{ab}$, there are hyperplanes $H_0$ and $H_1$ in $\Bbb P^N_k$,
such that $C\cap H_0\cap H_1=\emptyset$ and
$$\lambda_{H_0,v}(y) + \lambda_{H_1,v}(y)
  \ge \beta\alpha\lambda_{ij,v}(y) + O(1)\tag\02.5.3.2$$
for all $v\in S$ and all but finitely many $y\in C(k)$.
\endroster
\endit

\demo{Proof}  Let $\sigma_i$ and $\sigma_j$ be the canonical sections of
$\Cal O(A_i)$ and $\Cal O(A_j)$, respectively.  Then the linear system
$$\sigma_i^{\beta\alpha}\cdot\Gamma(Y,\beta(M-\alpha A_i))
  + \sigma_j^{\beta\alpha}\cdot\Gamma(Y,\beta(M-\alpha A_j))$$
has base locus $\Supp A_i\cap\Supp A_j$, since the first summand has
base locus $\Supp A_i$ and the second has base locus $\Supp A_j$.
This intersection consists of the union of $Z_{ij}$ and finitely many closed
points.
Choose an element of this linear system, sufficiently generic so that
it does not vanish identically on any irreducible component of $Z_{\ell m}$,
and let $H_0$ be the associated hyperplane in $\Bbb P^N_k$.  Fix a
Weil function $\lambda_{H_0}$ associated to $H_0$; since
$H_0-\beta\alpha Z_{ij}$ is an effective divisor, we have
$$\lambda_{H_0,v}(y) \ge \beta\alpha\lambda_{ij,v}(y) + O(1)
  \tag\02.5.3.3$$
for all $y\in Y(k)\setminus H_0$.

Similarly let $\sigma_\ell$ and $\sigma_m$ be the canonical sections of
$\Cal O(A_\ell)$ and $\Cal O(A_m)$, and let $H_1$ be the hyperplane associated
to an element of
$$\sigma_\ell^{\beta\alpha}\cdot\Gamma(Y,\beta(M-\alpha A_\ell))
  + \sigma_m^{\beta\alpha}\cdot\Gamma(Y,\beta(M-\alpha A_m))\;,$$
chosen sufficiently generically such that $H_1$ does not contain any
irreducible component of $H_0\cap Y$.  We also have
$$\lambda_{H_1,v}(y) \ge \beta\alpha\lambda_{\ell m,v}(y) + O(1)
  \tag\02.5.3.4$$
for all $y\in Y(k)\setminus H_1$.

By construction, $Y\cap H_0\cap H_1$ is a finite union of closed points,
so we can let $H_2$ be a hyperplane that avoids those points to ensure that
$Y\cap H_0\cap H_1\cap H_2=\emptyset$.  Since $\lambda_{H_2,v}\ge O(1)$,
(\02.5.3.1) follows from (\02.5.3.3) and (\02.5.3.4).  This proves (a).

For part (b), let $\sigma_i$ be as above, and let $H_0$ be the hyperplane
associated to an element of
$\sigma_i^{\beta\alpha}\cdot\Gamma(Y,\beta(M-\alpha A_i))$, chosen
generically so that $H_0$ does not contain $C$.  Let $H_1$ be a
hyperplane in $\Bbb P^N_k$, chosen so that $C\cap H_0\cap H_1=\emptyset$.
The choice of $H_0$ implies that
$$\lambda_{H_0,v}(y) \ge \lambda_{ij,v}(y) + O(1)$$
for all but finitely many $y\in C(k)$, so (\02.5.3.2) holds.\qed
\enddemo

\lemma{\02.5.4}  Let $i,j,\ell$ be distinct indices.  Then:
\roster
\myitem a.  there exist hyperplanes $H_0$, $H_1$, and $H_2$ in $\Bbb P^N_k$,
such that
$$Y\cap H_0\cap H_1\cap H_2=\emptyset$$
and
$$\lambda_{H_0,v}(y) + \lambda_{H_1,v}(y) + \lambda_{H_2,v}(y)
  \ge \beta\alpha(\lambda_{ij,v}(y) + \lambda_{i\ell,v}(y)
    + \lambda_{j\ell,v}(y)) + O(1)$$
for all $v\in S$ and all $y\in Y(k)$ outside of $H_0\cup H_1\cup H_2$; and
\myitem b.  given any integral curve $C\subseteq Y$ not contained in any of
the $Z_{ab}$, there are hyperplanes $H_0$ and $H_1$ in $\Bbb P^N_k$,
such that $C\cap H_0\cap H_1=\emptyset$ and
$$\lambda_{H_0,v}(y) + \lambda_{H_1,v}(y)
  \ge \beta\alpha(\lambda_{ij,v}(y) + \lambda_{i\ell,v}(y)) + O(1)
  \tag\02.5.4.1$$
for all $v\in S$ and all but finitely many $y\in C(k)$.
\endroster
\endit

\demo{Proof}  Let $\sigma_i$ and $\sigma_j$ be as in the preceding proof.
Choose a section of the linear system
$$\sigma_i^{\beta\alpha}\cdot\Gamma(Y,\beta(M-\alpha A_i))
  + \sigma_j^{\beta\alpha}\cdot\Gamma(Y,\beta(M-\alpha A_j))\;,$$
and let $H_0$ be the associated hyperplane.  We may assume that the choice
is sufficiently generic so that $H_0$ does not contain any irreducible
component of $A_\ell$.  We have
$$\lambda_{H_0,v}(y) \ge \beta\alpha\lambda_{ij,v}(y) + O(1)$$
for all $y\in Y(k)\setminus H_0$.

Next let $\sigma_\ell$ be the canonical section of $\Cal O(A_\ell)$,
and let $H_1$ be the hyperplane associated to a section of
$$\sigma_\ell^{\beta\alpha}\cdot\Gamma(Y,\beta(M-\alpha A_\ell))\;.$$
We may assume that $H_1$ does not contain any irreducible component
of $Y\cap H_0$.  We have
$$\lambda_{H_1,v}(y)
  \ge \beta\alpha(\lambda_{i\ell,v}(y) + \lambda_{j\ell,v}(y)) + O(1)\;.$$

Again, $Y\cap H_0\cap H_1$ consists of finitely many points, and we choose
$H_2$ to be any hyperplane not meeting any of these points.  Part (a)
then concludes as in the previous lemma.

For part (b), let $H_0$ and $H_1$ be the hyperplanes associated to suitably
chosen sections of
$\sigma_i^{\beta\alpha}\cdot\Gamma(Y,\beta(M-\alpha A_i))$ and
$\Gamma(Y,\beta M)$, respectively.  As in the previous lemma, we then have
$C\cap H_0\cap H_1=\emptyset$; moreover the choice of $H_0$ implies
$$\lambda_{H_0,v}(y)
  \ge \beta\alpha(\lambda_{ij,v}(y) + \lambda_{i\ell,v}(y)) + O(1)$$
for all but finitely many $y\in C(k)$, giving (\02.5.4.1).\qed
\enddemo

Now consider part (a) of the proposition.

The conditions in (\02.4) on the intersections of the divisors $Z_{ij}$
imply that there is a constant $C_v$ such that, for each $y\in Y(k)$
not in $\bigcup\Supp Z_{ij}$, one of the following conditions holds.
\roster
\myitem i.  $\lambda_{ij,v}(y)\le C_v$ for all $i$ and $j$;
\myitem ii.  there are indices $i$ and $j$ such that $\lambda_{ij,v}(y)>C_v$
but $\lambda_{ab,v}(y)\le C_v$ in all other cases;
\myitem iii.  there are distinct indices $i,j,\ell,m$ such that
$\lambda_{ij,v}(y)>C_v$ and $\lambda_{\ell m,v}(y)>C_v$\snug\hbox{,\kern-1pt}
but $\lambda_{ab,v}(y)\le C_v$ in all other cases; or
\myitem iv.  there are indices $i,j,\ell$ such that
$\max\{\lambda_{ij,v}(y), \lambda_{i\ell,v}(y), \lambda_{j\ell,v}(y)\}>C_v$,
but $\lambda_{ab,v}(y)\le C_v$ if $\{a,b\}\nsubseteq\{i,j,\ell\}$.
\endroster

In cases (iii) and (iv), (\02.5.1) follows from Lemmas \02.5.3a and \02.5.4a,
respectively.  In case (i) there is nothing to prove.  Case (ii) follows as
a special case of Lemma \02.5.3a, since $n\ge4$.

Since the indices in Lemmas \02.5.3 and \02.5.4 have only finitely many
possibilities, the set of hyperplanes that occur can be assumed to be finite.
This proves (a).

For part (b), we have cases (i)--(iv) as before.  Cases (ii) and (iii)
follow from Lemma \02.5.3b, where we may assume without loss of generality
that $\lambda_{ij,v}(y)\ge\lambda_{\ell m,v}(y)$ to obtain (\02.5.2) from
(\02.5.3.2).  Similarly, case (iv) follows from Lemma \02.5.4b after a suitable
permutation of the indices, and case (i) is again trivial.  The set of
hyperplanes can again be taken to be finite, for the same reason.\qed
\enddemo

The main theorem of this paper can now be stated and proved.

\thm{\02.6}  Let $k$ be a number field, let $S$ be a finite set of places
of $k$, and let $Y$, $n$, $\{Z_{ij}\}_{i<j}$, $\{A_i\}_i$, $M$, and $\alpha$
be as in (\02.4).  Then:
\roster
\myitem a.  if $\alpha>6$ then no set of $\Cal O_{k,S}$\snug-integral points
on $Y\setminus\bigcup Z_{ij}$ is Zariski-dense, and
\myitem b.  if $\alpha>8$ then every set of $\Cal O_{k,S}$\snug-integral points
on $Y\setminus\bigcup Z_{ij}$ is finite.
\endroster
\endit

\demo{Proof}  Let $\beta$ be as in Proposition \02.5, and again regard $Y$
as a subvariety of $\Bbb P^N_k$ as in that proposition.  For points $y\in Y(k)$
let $h(y)$ denote the height of $y$ as a point in $P^N_k$.  Fix Weil functions
$\lambda_{ij}$ for each $Z_{ij}$.

Assume that $\alpha>6$.  Note that then $n\ge4$, since
$$\sum_{i=1}^n(M-\alpha A_i) = nM-\alpha\sum A_i = (n-\alpha)M$$
is ample, hence $n>\alpha>6$.

Assume, by way of contradiction, that some set of $\Cal O_{k,S}$\snug-integral
points $y$ on $Y\setminus\bigcup Z_{ij}$ is Zariski-dense.  Then, for
these integral points,
$$\sum_{v\in S}\sum_{i<j} \lambda_{ij,v}(y) = \frac1{2\beta}h(y) + O(1)\;,
  \tag\02.6.1$$
with the constant in $O(1)$ independent of $y$.

Let $H_1,\dots,H_q$ be as in Proposition \02.5a, and let $\lambda_{H_j}$ be
associated Weil functions for them.  Without loss of generality we may assume
that none of the integral points $y$ lie on any of these hyperplanes.
Then, by Proposition \02.5a, for each such $y$ and each $v\in S$ the inequality
$$\max_{j\in\Cal J} \sum_{j\in J} \lambda_{H_j,v}(y)
  \ge \beta\alpha\sum_{i<j} \lambda_{ij,v}(y) + O(1)$$
holds, where $\Cal J$ is as in Proposition \02.4a and where the constant
in $O(1)$ independent of $y$.

Combining this with (\02.6.1) then gives
$$\sum_{v\in S} \max_{j\in\Cal J} \sum_{j\in J} \lambda_{H_j,v}(y)
  \ge \frac\alpha2 h(y) + O(1)\;,$$
contradicting (\02.3).  This proves part (a).

Now consider part (b).  By part (a), it will suffice to assume that some
integral curve $C$ on $Y$ contains infinitely many integral points, and derive
a contradiction.  Such a $C$ cannot be contained in the support of any $Z_{ij}$,
so Proposition \02.5b applies.

Therefore, there are finitely many hyperplanes $H_1,\dots,H_q$ in $\Bbb P^N_k$
such that, for each integral point $y$ and each $v\in S$ the inequality
$$\max_{j\in\Cal J} \sum_{j\in J} \lambda_{H_j,v}(y)
  \ge \frac{\beta\alpha}{2} \sum_{i<j} \lambda_{ij,v}(y) + O(1)$$
holds, where $\Cal J$ is as in Proposition \02.4b and with $O(1)$ independent
of $y$.  Again, combining this with (\02.6.1) gives
$$\sum_{v\in S} \max_{j\in\Cal J} \sum_{j\in J} \lambda_{H_j,v}(y)
  \ge \frac\alpha4 h(y) + O(1)\;.$$
This again contradicts (\02.3), since $\alpha>8$.\qed
\enddemo

\beginsection{\03}{The Nevanlinna Case}

The Nevanlinna counterpart to Theorem \02.6 can be proved using
substantially the same method, using the Nevanlinna counterpart to
the Evertse-Ferretti theorem (Theorem \02.1) due to M. Ru \cite{R}.
The statement here is a slight variation on that theorem; see
\cite{V~2, Thm.~21.7}.

\thm{\03.1} (Ru)  Let $X$ be a closed subvariety of $\Bbb P^N_{\Bbb C}$
of dimension $n\ge1$, let $D_0,\dots,D_q$ be hypersurfaces in
$\Bbb P^N_{\Bbb C}$, let $\lambda_{D_0},\dots,\lambda_{D_q}$ be corresponding
Weil functions, and let $\epsilon>0$.  Let $\Cal J$ be the set of
all $(n+1)$\snug-element subsets $J$ of $\{1,\dots,q\}$ for which
$$X\cap\bigcap_{j\in J} D_j = \emptyset\;,$$
and assume that $\Cal J$ is nonempty.  Finally, let $f\:\Bbb C\to X$ be
a holomorphic function with Zariski-dense image.  Then
$$\int_0^{2\pi} \max_{J\in\Cal J} \sum_{j\in J}
    \frac{\lambda_{D_j}(f(re^{\sqrt{-1}\theta}))}{\deg D_j}\,
    \frac{d\theta}{2\pi}
  \le_\exc (n+1+\epsilon)T_f(r) + O(1)\;.\tag\03.1.1$$
\endit

Here the notation $\le_\exc$ means that the inequality holds for all $r>0$
except for a set of finite Lebesgue measure.

Again, we only need the case in which the $f_j$ are all linear, associated
to hyperplanes $H_j$ with Weil functions $\lambda_{H_j}$.  In this case,
(\03.1.1) reduces to
$$\int_0^{2\pi} \max_{J\in\Cal J}\sum_{j\in J}
     \lambda_{H_j}(f(re^{\sqrt{-1}\theta}))\,\frac{d\theta}{2\pi}
  \le_\exc (n+1+\epsilon)T_f(r) + O(1)\;.\tag\03.2$$

The following theorem corresponds to Theorem \02.6.

\thm{\03.3}  Let $k=\Bbb C$ and let $Y$, $n$, $\{Z_{ij}\}_{i<j}$,
$\{A_i\}_i$, $M$, and $\alpha$ be as in (\02.4).
Let $f\:\Bbb C\to Y\setminus\bigcup\Supp Z_{ij}$ be a holomorphic curve.  Then:
\roster
\myitem a.  if $\alpha>6$ then the image of $f$ is not Zariski-dense.
\myitem b.  if $\alpha>8$ then $f$ must be constant.
\endroster
\endit

\demo{Proof}  Let $\beta$ and $Y\hookrightarrow\Bbb P^N_{\Bbb C}$ be as in
the proof of Theorem \02.6, and let $T_f(r)$ be defined via this embedding.
For each $i<j$ let $\lambda_{ij}$ be a Weil function for $Z_{ij}$.  Then
$$\int_0^{2\pi} \sum_{i<j} \lambda_{ij}(f(re^{\sqrt{-1}\theta}))\,
    \frac{d\theta}{2\pi}
  = \frac1{2\beta}T_f(r) + O(1)\;.\tag\03.3.1$$
As before, the assumptions imply that $n\ge4$; then Proposition \02.5a applies,
giving hyperplanes $H_1,\dots,H_q$ in $\Bbb P^N_{\Bbb C}$.  Let $\Cal J$ be
the set of $3$\snug-element subsets $J$ of $\{1,\dots,q\}$ for which
$X\cap\bigcap_{j\in J}H_j=\emptyset$, and note that $\Cal J\ne\emptyset$.
For all $z\in\Bbb C$ outside of a discrete subset, we have
$f(z)\notin\bigcup H_j$, and for those $z$ we have
$$\max_{j\in\Cal J}\sum_{j\in J}\lambda_{H_j}(f(z))
  \ge \beta\alpha\sum_{i<j} \lambda_{ij}(f(z)) + O(1)\;.$$
Combining this with (\03.3.1) gives
$$\int_0^{2\pi} \max_{J\in\Cal J}\sum_{j\in J}
     \lambda_{H_j}(f(re^{\sqrt{-1}\theta}))\,\frac{d\theta}{2\pi}
  \ge \frac\alpha 2 T_f(r) + O(1)\;.$$
If the image of $f$ is Zariski-dense, then this contradicts (\03.2)
since $\alpha>6$.  This proves (a).

The proof of (b) is similar to the proof of Theorem \02.6b, with the same
types of changes as for part (a).  It is left to the reader.\qed
\enddemo

\enddocument